\documentclass[12pt]{article}
\usepackage{amsfonts}
\usepackage{amssymb}
\usepackage{amsmath}
\usepackage{eufrak}
\usepackage[dvips]{graphicx}
\usepackage{latexsym}
\newcommand{\ind}{\makebox[1em]{\raisebox{-.5ex}[0ex][0ex]{\makebox[0em]%
{$\smile$}}\raisebox{.4ex}[0ex][0ex]{\makebox[-.02em]{$|$}}}}

\newcommand{\dep}{\makebox[1em]{\raisebox{.3ex}[0ex][0ex]%
{$\not$}\makebox[.7em]{\ind}}}

\newcommand{\nmdep}{\makebox[1em]{\raisebox{1.5ex}[0ex][0ex]{\makebox[0em]%
{$\scriptscriptstyle n\! m$}}\makebox[-1em]{$\dep$}}}
\newcommand{\nmind}{\makebox[1em]{\raisebox{1.5ex}[0ex][0ex]{\makebox[0em]%
{$\scriptscriptstyle n\! m$}}\makebox[-1em]{$\ind$}}}

\newcommand{\bi}{\begin{itemize}}
\newcommand{\ei}{\end{itemize}}

\newtheorem{theorem}{Theorem}[section]
\newtheorem{example}[theorem]{Example}

\newtheorem{fact}[theorem]{Fact}
\newtheorem{claim}{Claim}
\newtheorem{corollary}[theorem]{Corollary}
\newtheorem{proposition}[theorem]{Proposition}
\newtheorem{definition}[theorem]{Definition}
\newtheorem{remark}[theorem]{Remark}

\newtheorem{question}[theorem]{Question}

\newtheorem{problem}[theorem]{Problem}

\title{Topologies induced by group actions}
\author{Jan Dobrowolski\footnote{Research supported by NCN Grant no. 2012/05/N/ST1/02850}}
\date{}

\begin{document}
\maketitle
\begin{abstract} We introduce some canonical topologies induced by
actions of topological groups on groups and rings. For $H$ being a group [or a ring]
and $G$ a topological group acting on $H$ as automorphisms, we describe the
finest group [ring] topology on $H$ under which the action of $G$ on
$H$ is continuous. We also study the introduced topologies
in the context of Polish structures. In particular, we prove 
that there may be no Hausdorff topology on a group $H$ under which
a given action of a Polish group on $H$ is continuous.
\end{abstract}

\section{Introduction}
The main motivation for this paper is the following general problem. Suppose $G$ is a topological group
acting on $X$,
where $X$ is a set, possibly equipped with some algebraic structure preserved by the action of $G$.
When does there exist a "nice" topology on $X$, such that the action of $G$ on $X$ is continuous,
and the topology is compatible with the structure on $X$? Clearly, if there is such a topology which is
at least $T_1$, then, for every element $x\in X$, its stabilizer $G_x$ is closed in $G$. On the other hand,
if the latter is satisfied, then, by Remark \ref{tau} below, 
we can equip $X$ with a topology under which the action is continuous and
which inherits many properties of the given topology on $G$.

Now, suppose that $X$ is equipped with a group structure (preserved by the action of $G$).
Then, the topology $\tau$ defined below
usually fails to be a group topology. In Theorem \ref{main},
we give a
description of the finest group topology on $X$, under which 
the action of $G$ on $X$ is continuous. Using this description, we give an example of
an action of the polish group $Homeo([0,1])$ on a certain group $H$, such that there
is no Hausdorff group topology on $H$ under which the action is continuous.

Also, we give in Theorem \ref{ring_top} a description of the finest compatible topology in the 
case of $X$ being a ring.

By a topological group we will mean a group equipped with a topology, such that
the multiplication and the inversion are continuous functions (we do
not assume that the topology is Hausdorff).

\begin{remark}\label{tau}
Suppose $G$ is a topological group acting on a set $X$. 
Define ${\cal U}:=\{ U\cdot x: U \subseteq G, U$ is open, $ x \in X\}$. Then, we have:\\
(1) The family ${\cal U}$ is a basis of a topology on $X$. Denote this topology by $\tau$.\\
(2) The action of $G$ on $(X,\tau)$ is continuous.\\
(3) All $G$-orbits in $X$ are clopen in $\tau$; 
moreover, for every $x \in X$, $G/G_x \approx G\cdot x$ (where the orbit $G\cdot x$ is equipped with 
the topology $\tau$).\\
(4) $(\forall x \in X) (G_x$ is closed in $G)$ $\iff$ $\tau$ is $T_1$ $\iff$ $\tau$ 
is Hausdorff.
\end{remark}
{\em Proof.}\\
(1) Take any $x_1,x_2\in X$, open sets $U_1,U_2\subseteq G$ and a point $y\in U_1x_1\cap U_2x_2$.
Then, $y=u_1x_1=u_2x_2$ for some $u_1\in U_1$ and $u_2\in U_2$. For $i=1,2$, let $W_i$ be 
an open neighbourhood of $e$ in $G$ such that $W_iu_i\subseteq U_i$. Put $W=W_1\cap W_2$. 
Then, $Wy=Wu_ix_i\subseteq U_ix_i$ for $i=1,2$, so we are done.\\
(2) Take any $g \in G$, $x,y\in X$ and an open set $U\subseteq G$ such that $gx\in Uy$. Choose $u\in U$
so that $gx=uy$. Let $U_1$ be an open neighbourhood of $e$ in $G$ such that $U_1u\subseteq U$, and
let $V_1,V_2\subseteq G$ be open neighbourhoods of $e$ such that $V_1gV_2\subseteq U_1g$. Then,
$(V_1g)(V_2x)\subseteq U_1gx=U_1uy\subseteq Uy$, which yields the continuity of the action.\\
(3) By the definition of $\tau$, every $G$-orbit is open, so also clopen. 
It is straightforward to check that $f_x:G/G_x\to Gx$ given by $aG_x\mapsto ax$ is a homeomorphism.\\
(4) First condition implies third by (3) and the fact that the quotient of
any topologcial group (not necessarily Hausdorff) by a closed subgroup is Hausdorff. 
Third condition implies second
trivially, and second implies first by (2).
\hfill $\square$\\
We will denote the topology $\tau$ from the above remark by $\tau(X,G)$.

\section{Group actions}\label{groups}

In this section, we will describe the finest group [ring] topology on $H$ (where $H$ is 
a group or a ring, respectively) under which a given action of a topological group on $H$ (as
automorphisms) is continuous.

First, let us recall a recent result of Bergman, which we use in our construction.
Let $G$ be a group. We will denote the neutral element of $G$ by $e$, and 
for $S\subseteq G$ we will write $S^*=S\cup \{e\}\cup S^{-1}$.
Following the notation from \cite{83}, given any $\mathbb{Q}$-tuple
$(S_q)_{q\in \mathbb{Q}}$, put
$$U((S_q)_{q\in \mathbb{Q}})=\bigcup_{n<\omega}\bigcup_{q_1<\dots<q_n\in\mathbb{Q}}S_{q_1}^*\dots S_{q_n}^*.$$
Below, we will omit the symbol $\bigcup_{n<\omega}$ in similar expressions.
For a family $F$ of subsets of $G$, $F^G$ will denote the collection of all subsets of $G$ of 
the form $\bigcup_{g\in G} gS_gg^{-1}$, for $G$-tuples $(S_g)_{g\in G}$ of members of $F$.
We say that a filter on $G$ converges to $e$ in a given topology, if every neighbourhood
of $e$ contains a member of the filter.
By Lemma 14 and Proposition 15 from \cite{83}, we have:
\begin{fact}\label{bergman}
Let $F$ be a downward directed family of nonempty subsets of $G$. Then the sets
$U((S_q)_{q\in \mathbb{Q}})$, where $(S_q)_{q\in \mathbb{Q}}$ ranges over all $\mathbb{Q}$-tuples
of members of $F^G$, form a basis of open neighbourhoods of $e$ in a group topology $\mathcal{T}_F$, 
which is the finest group topology on $G$ under which $F$ converges to $e$.
\end{fact}
When $\rho$ is a topology on $G$, we will denote by $\rho ^*$ the topology $\mathcal{T}_F$, 
where $F$ consists of
$\rho$-open neighbourhoods of $e$. In particular, if $\rho$ is a group topology, then $\rho^*=\rho$.

Let $G$ be a topological group equipped with a topology $\sigma$.

\begin{theorem}\label{main}
Suppose $G$ acts on a group $H$ as automorphisms. 
We identify $G$ and $H$ with $\{e\}\times G<H\rtimes G$ and $H\times\{e\}<H\rtimes G$, respectively.
Put $T=(D\times \sigma)^*$, where $D$ is
the discrete topology on $H$. Denote by $T_H$ and $T_G$ the topologies
induced by $T$ on the subgroups identified with $H$ and $G$, respectively.
Then:\\
(1) $T_G=\sigma$.\\
(2) $T_H$ is a group topology on $H$ under which the action of $G$ on $H$ is continuous.\\
(3) If $\rho$ is another group topology on $H$ under which the action of $G$ on $H$ is continuous,
then $T_H$ is finer than $\rho$.\\
We will denote the topology $T_H$ by $T(H,G)$.
\end{theorem}
{\em Proof.}\\
(1) We will show that $T_G =\sigma^*$ (using the description of $T=(D\times \sigma)^*$ and $\sigma^*$
given by Fact \ref{bergman}), which suffices since $\sigma=\sigma^*$.
It is easy to see that $T_G$ is a group topology on $G$. Hence, it is enough to show that the neighbourhoods
of $e$ in $T_G$ are the same as in $\sigma^*$.

Take a $T_G$-open neighbourhood of $e$ of the form 
$$V=G\cap\bigcup_{q_1<\dots<q_n\in\mathbb{Q}}\Big( \bigcup_{(h,g)\in H\rtimes G}(\{e\}\times U_{(h,g)}^{q_1})^{(h,g)}\Big)
\dots \Big(\bigcup_{(h,g)\in H\rtimes G}(\{e\}\times U_{(h,g)}^{q_n})^{(h,g)}\Big)=$$\\
$$=G\cap\bigcup_{q_1<\dots <q_n\in\mathbb{Q}}\Big(\bigcup_{(h,g)\in H\rtimes G}
\{(h(u^gh^{-1}),u^g):u\in U_{(h,g)}^{q_1}\}\Big)\dots$$\\
$$\dots\Big(\bigcup_{(h,g)\in H\rtimes G}\{(h(u^gh^{-1}),u^g):u\in U_{(h,g)}^{q_n}\}\Big),$$ 
where each $U_{(h,g)}^q$ is a symmetric $\sigma$-open neighbourhood of $e$ in $G$.
Then,
$$\bigcup_{q_1<\dots <q_n}\Big( \bigcup_{g\in G}gU_{(e,g)}^{q_1}g^{-1}\Big)\dots 
\Big(\bigcup_{g\in G}gU_{(e,g)}^{q_n}g^{-1}\Big)$$ is a 
$\sigma^{*}$-open neighbourhood of $e$ 
contained in $V$.

Conversly, take a $\sigma^{*}$-open neighbourhood of $e$ of the form 
$$W=\bigcup_{q_1<\dots <q_n}\Big( \bigcup_{g\in G}gU_g^{q_1}g^{-1}\Big)\dots\Big(\bigcup_{g\in G}gU_g^{q_n}g^{-1}\Big),$$
where each $U_g^q$ is a symmetric neighbourhood of $e$ in $G$.
For any $(h,g)\in G$ and $q\in \mathbb{Q}$, find a $\sigma$-open, symmetric neighbourhood $U_{(h,g)}^q$ of $e$, whose 
conjugate by $g$ is contained in $U_g$ (it can be chosen independently from $h$). Then, 
$$G\cap\bigcup_{q_1<\dots<q_n\in\mathbb{Q}}\Big(\bigcup_{(h,g)\in H\rtimes G}(\{e\}\times U_{(h,g)}^{q_1})^{(h,g)}\Big)
\dots \Big(\bigcup_{(h,g)\in H\rtimes G}(\{e\}\times U_{(h,g)}^{q_n})^{(h,g)}\Big)$$
is a $T_G$-open neighbourhood of $e$ contained in $W$.\\
(2) $T_H$ is a group topology since $H$ is a subgroup of $H\rtimes G$, and
$T$ is a group topology on $H\rtimes G$ by Fact \ref{bergman}. For the continuity of the action,
take any
$g\in G$, $h\in H$ and a $T$-open set $U$, such that $gh\in U\cap H$. This means that, in $H\rtimes G$,
$(e,g)(h,e)(e,g)^{-1}\in U$, so we can choose $T$-open sets $U_1$ and $U_2$, such that $(e,g)\in U_1$,
$(h,e)\in U_2$ and $U_1U_2U_1^{-1}\subseteq U$.
Then, for any $g_1\in U_1\cap G$ and $h_1$ in $U_2\cap H$, we have that $(e,g_1)(h_1,e)(e,g_1^{-1})=(g_1h_1,e)$ belongs 
to $U$, so $g_1h_1\in U\cap H$. This proves the continuity of the action.\\
(3) Suppose $\rho$ is a group topology on $H$ under which the action of $G$ on $H$ is continuous. Then,
 the product topology $\rho\times\sigma$ is a group topology on $H\rtimes G$, which is coarser
 than $D\times\sigma$,
 so, by the choice of $T$, we have that $T$ is finer than $\rho\times\sigma$. In particular, $T_H$ is finer
 than $\rho$.
\hfill $\square$\\

Using the above theorem we obtain an explicit formula describing the
topology $T(H,G)$:
\begin{corollary}\label{description}
With the notation from the above theorem, $T(H,G)$ has a basis of open neighbourhoods of $e$ consisting of the sets:
$$\bigcup_{q_1<\dots<q_n\in\mathbb{Q}}\{h_1(u_1h_1^{-1})u_1(h_2(u_2h_2^{-1}))u_1u_2(h_3(u_3h_3^{-1}))
\dots u_1u_2\dots u_{n-1}(h_n(u_nh_n^{-1})):$$\\
$$h_i\in H,u_i\in U_{h_i}^{q_i}, u_1\dots u_n=e\},$$
where $(U_h^q)_{h\in H,q\in \mathbb{Q}}$ range over all 
$H\times \mathbb{Q}$-tuples of $\sigma$-open symmetric neighbourhoods of $e$ in $G$.
\end{corollary}
{\em Proof.}
By the description of the topology $T_H$ given in Fact \ref{bergman}, we get that it
has a basis of open neighbourhoods of $e$ consisting of the sets:
$$\bigcup_{q_1<\dots<q_n\in\mathbb{Q}}\{h_1(v_1^{g_1}h_1^{-1}) v_1^{g_1}(h_2(v_2^{g_2}h_2^{-1}))
\dots (v_1^{g_1}\dots 
v_{n-1}^{g_n})(h_n(v_n^{g_n}h_n^{-1})):
h_i\in H,$$ \\ $$g_i\in G, v_i\in U_{(h_i,g_i)}^{q_i}, v_1^{g_1}\dots v_n^{g_n}=e\}=$$\\
$$=\bigcup_{q_1<\dots<q_n\in\mathbb{Q}}\{h_1(u_1h_1^{-1}) u_1(h_2(u_2h_2^{-1}))\dots (u_1\dots 
u_{n-1})(h_n(u_nh_n^{-1})):
h_i\in H,$$ \\ $$g_i\in G, u_i\in (U_{(h_i,g_i)}^{q_i})^{g_i}, u_1\dots u_n=e\},$$
where $(U_{(h,g)}^q)_{(h,g)\in H\rtimes G,q\in \mathbb{Q}}$ range over all 
$(H\rtimes G)\times \mathbb{Q}$-tuples of $\sigma$-open symmetric neighbourhoods of 
$e$ in $G$.
Since the tuples $(U_{(h,g)}^{q})^{g}_{(h,g)\in H\rtimes G,q\in \mathbb{Q}}$ range over the 
same set, we can omit the conjugations in the formula.
Now, if we replace each $U_{(h,g)}^q$ by $U_{(h,e)}^q$ in a tuple 
$(U_{(h,g)}^q)_{(h,g)\in H\rtimes G,q\in \mathbb{Q}}$,
then the corresponding neighbourhood of $e\in H$ will be contained in
the original one. Thus, we obtain
the same topology when we restrict ourselves to tuples in 
which $U_{(h,g)}^q=U_h^q$ does not depend on $g$.
This gives the conclusion.
\hfill $\square$\\

In Section \ref{topologies}, we will use the description of $T(H,G)$ that we have obtained
to prove the absence of a compatible Hausdorff topology
for some classes of Polish group structures (see Proposition \ref{haus}).

Let us keep the notation from above and define $\lambda(H,G)$ to be the 
topology on $H$ in which a set $U$ is open if
for each $h_1,h_2\in H$, the sets $h_1Uh_2,h_1U^{-1}h_2$ are open in 
the topology $\tau(H,G)$ (defined after Remark \ref{tau}).
It is easy to see that if we equip $H$ with $\lambda(H,G)$, 
then the action of $G$ on $H$ is separately continuous,
the inversion on $H$ is continuous and the multiplication on $H$ is 
separately continuous. Moreover, $\lambda(H,G)$ is the finest topology on $H$ 
with these properties. Indeed, let $\xi$ be any other such topology. 
Take any $\xi$-open set $V$. Then, for any $h_1,h_2\in H$,
$h_1Uh_2,h_1U^{-1}h_2$ are $\xi$-open,
so also $\tau$-open. Hence, $V$ is $\lambda(H,G)$-open.

\begin{remark}\label{other_top}
In Theorem \ref{main}, we can replace the discrete topology $D$ by any topology on $H$ which is
finer than all group topologies under which the action of $G$ on $H$ is continuous. Examples 
of such topologies are $\tau(H,G)$ 
and $\lambda(H,G)$. 
However, the simplest description of $T(H,G)$ we 
obtain starting from the
discrete topology on $H$.
\end{remark}
Let us formulate a remark about the topology $\lambda(H,G)$ defined above.

\begin{remark}
If the topology of $G$ and $\lambda(H,G)$ are metrizable and Baire, then\\
$\lambda(H,G)=T(H,G)$
\end{remark}
{\em Proof.}
Let us equip $H$ with the topology $\lambda(H,G)$.
Since the multiplication on $H$ is separately conitnuous, the inversion on $H$ is continuous and
the action of $G$ on $H$ is separately continuous, we get by
Theorem 9.14 from \cite{87} that
$H$ is a topological
group with the 
topology $\lambda(H,G)$, and the action of $G$ on $H$ is continuous. So, $\lambda(H,G)$ is coarser
than $T(H,G)$. But $\lambda(H,G)$ is always finer than $T(H,G)$ (see the discussion
preceeding Remark \ref{other_top}), so these topologies are equal.
\hfill $\square$\\
The above remark can by illustrated by the following example.
\begin{example}
Let $G=S_\omega$ be the group of all permutations of $\omega$, considered with the product topology (which
is Polish, so, in particular, metrizable and Baire). Consider the action of $G$ on $H:=2^{\omega}$,
given by $g\cdot h=h\circ g^{-1}$. Then, $\lambda(H,G)$ is the product topology on $2^{\omega}$,
so it coincides with $T(H,G)$.
\end{example}
{\em Proof.}
Clearly $\lambda(H,G)$ is finer than the product topology. For the converse, let $U$
be any $\lambda(H,G)$-open neighbourhood of $0\in H$. We will show that it contains
an  open neighbourhood of $0\in H$ in the sense of the product topology.
Let $\omega=A\dot{\cup} B\dot{\cup} C$ be a partition of $\omega$ into three infinite sets.
For $i,j,k \in \{0,1\}$ define $h_{i,j,k}\in H$ to be equal to $i$ on $A$, equal to $j$ on $B$,
and equal to $k$ on $C$. 
For $i,j,k \in \{0,1\}$, $h_{i,j,k}+U$ contains a $\tau(H,G)$-open neighbourhood of $h_{i,j,k}$
of the form $[\alpha_{i,j,k}]\cdot h_{i,j,k}$, where $\alpha_{i,j,k}:\omega\to \omega$ is a partial
function with a finite domain, and $[\alpha]=\{\eta\in S_{\omega}:\alpha\subseteq \eta\}$.
We finish by the following claim:
\begin{claim}
Put $I=\bigcup_{i,j,k\in \{0,1\}} (dom(\alpha_{i,j,k})\cup rng(\alpha_{i,j,k}))$. Then,
$\{x\in H:x_{|I}=0\}\subseteq U.$
\end{claim}\label{cl}
{\em Proof of Claim \ref{cl}.}
Take any $x\in H$ such that $x_{|I}=0$.
Notice that we can choose $i,j,k\in \{0,1\}$ so that $(h_{i,j,k}+x)^{-1}[\{0\}]$ and 
$(h_{i,j,k}+x)^{-1}[\{1\}]$ are both infinite, and $i,j,k$ are not all equal.
Then, since $h_{i,j,k}^{-1}[\{0\}]$ and $h_{i,j,k}^{-1}[\{1\}]$ are also both infinite, and 
$h_{i,j,k}+x$ agrees with $h_{i,j,k}$ on $I$, we can find a permutation $\eta\in [\alpha_{i,j,k}]$, such that 
$\eta\cdot h_{i,j,k}=h_{i,j,k}+x$. Thus, $h_{i,j,k}+x\in h_{i,j,k}+U$, so $x\in U$.
\hfill $\square$\\

Now, we aim towards a description of the finest 
ring topology on $R$ under which a given action of a topological group on $R$
is continuous.
First, we give a variant of Fact \ref{bergman}, in which we are interested 
in semigroup topologies (i.e. topologies under which the multiplication is 
continuous) rather than group topologies on $G$, but still we assume
that $G$ is a group. For a subset $S$ of $G$, we will write $S^{\#}=S\cup \{e\}$.
We define
$$U'((S_q)_{q\in \mathbb{Q}})=\bigcup_{q_1<\dots<q_n\in\mathbb{Q}}S_{q_1}^{\#}\dots S_{q_n}^{\#}.$$
Then, by a straightforward modification (which is just replacing
expressions of the form $S^{*}$ by $S^{\#}$) of the proof of Lemma 14 and Proposition 15 
from \cite{83}, we obtain:
\begin{fact}\label{bergman'}
Let $F$ be a downward directed family of nonempty subsets of $G$. Then, the sets
$U'((S_q)_{q\in \mathbb{Q}})$, where $(S_q)_{q\in \mathbb{Q}}$ ranges over all $\mathbb{Q}$-tuples
of members of $F^G$, form a basis of open neighbourhoods of $e$ in a semigroup topology $\mathcal{T'}_F$, 
which is the finest semigroup topology on $G$ under which $F$ converges to $e$.
\end{fact}

When $\rho$ is a topology on $H$, we will denote by $\rho ^{\#}$ the topology $\mathcal{T}_F'$, 
where $F$ consist of
$\rho$-open neighbourhoods of $e$. 

Let $\sigma$ be a fixed group topology on a group $G$.
Repeating the proof of Theorem \ref{main},
we obtain:
\begin{proposition}\label{main'}
Suppose $G$ acts on a group $H$ as automorphisms. 
We identify $G$ and $H$ with $\{e\}\times G<H\rtimes G$ and $H\times\{e\}<H\rtimes G$, respectively.
Put $T'=(D\times \sigma)^{\#}$, where $D$ is
the discrete topology on $H$. Denote by $T_H'$ and $T_G'$ the topologies
induced by $T'$ on the subgroups identified with $H$ and $G$, respectively.
Then, $T_G'=\sigma$ and $T_H'$ is the finest semigroup topology on $H$ under which the action of $G$ on $H$ is continuous.
We will denote it by $T'(H,G)$.
\end{proposition}

Now, we are in a position to give a description of the finest topology in the ring case.
\setcounter{claim}{0}
\begin{theorem}\label{ring_top}
Suppose $G$ is a group equipped with a group topology $\sigma$, acting as 
automorphisms on a ring $R$.
Put $R_1=R\times \mathbb{Z}$, and define $+$ and $\cdot$ 
on $R_1$ by $(a,k)+(b,l)=(a+b,k+l)$ and
$(a,k)\cdot(b,l)=(ab+l\times a+k\times b,k\cdot l)$.
Clearly, $G$ acts on $R_1$ as automorphisms by $g(a,k):=(g(a),k)$.
Consider the induced action of $G$ on $GL_3(R_1)$.
We identify $R$ with a subset $\{\left( \begin{array}{ccc}
1 & 0 & x \\
0 & 1 & 0 \\
0 & 0 & 1 \end{array} \right):x\in R\}$ of $GL_3(R_1)$.
Denote by $T^r(R,G)$ the topology induced on $R$ by $T'(GL_3(R_1),G)$.
Then, $T^r(R,G)$ is the finest ring topology on $R$ such that the action
of $G$ on $R$ is continuous.

\end{theorem}
{\em Proof.}
Let $T_1$ be the topology induced on $R_1$ by $T'(GL_3(R_1),G)$ (we identify
$R_1$ with a subset of $GL_3(R_1)$ in the same manner as we do with $R$).
\begin{claim}
$T_1$ is the finest ring topology on $R_1$ under which the action of $G$ on $R_1$ is 
continuous.
\end{claim}
First, suppose the claim is proved and let us see that the conclusion of the theorem
follows.

By the claim, $T^r(R,G)$ is a ring topology on $R$, and the action of $G$ on $R$ 
equipped with $T^r(R,G)$ is continuous (as the restriction of the action on $R_1$).

Suppose $\rho$ is another topology on $R$ such that the action of $G$ on $R$ is 
continuous. Consider $R_1$ equipped with the product of $\rho$ and 
the discrete topology $E$ on $\mathbb{Z}$. Then, the  action $G$ on 
$R_1$ is also continuous and $R_1$ is a topological ring, so, by the claim, 
$T_1$ is finer than $\rho\times E$. Hence, $T^r(R,G)$ (which is equal to
the topology induced on $R$ by $T_1$) is finer than $\rho$, so we are done.\\
{\em Proof of Claim 1.}
First, we will check that $T_1$ is finer than every ring topology on $R_1$  under which the action
of $G$ on $R_1$ is continuous. Let $\chi$ be any such topology.
Let us equip $GL_3(R_1)$ with the topology $Z$
induced from the product topology $\chi^9$ on $R_1^9$. Then, $GL_3(R_1)$ becomes
a topological semigroup, and the action of $G$ on it is continuous.
Thus, $T'(GL_3(R_1),G)$ is finer that $Z$, so $T_1$ is finer than the 
topology induced by $Z$ on $R_1$, i.e. $T_1$ is finer than $\chi$.

Now, consider $R_1$ equipped with the topology $T_1$. The action of $G$ on $R_1$
is continuous (as a restriction of the action on $GL_3(R_1)$) and the addition in $R$
is continuous (as a restriction of the multiplication in $GL_3(R_1)$).
Moreover, the additive inversion in $R$ is continuous, as it is given by the map

$$\left( \begin{array}{ccc}
1 & 0 & x \\
0 & 1 & 0 \\
0 & 0 & 1 \end{array} \right)
\mapsto
\left( \begin{array}{ccc}
1 & 0 & -x\\
0 & 1 & 0 \\
0 & 0 & 1 \end{array} \right)=
\left( \begin{array}{ccc}
-1 & 0 & 0 \\
0 & 1 & 0 \\
0 & 0 & 1 \end{array} \right)
\left( \begin{array}{ccc}
1 & 0 & x \\
0 & 1 & 0 \\
0 & 0 & 1 \end{array} \right)
\left( \begin{array}{ccc}
-1 & 0 & 0 \\
0 & 1 & 0 \\
0 & 0 & 1 \end{array} \right),$$
which is continuous with respect to $T'(GL_3(R_1),G)$.

 It remains 
to show that the multiplication on $R_1$ is continuous. So, we will be done if 
we show that the map $$\Bigg(\left( \begin{array}{ccc}
1 & 0 & x \\
0 & 1 & 0 \\
0 & 0 & 1 \end{array} \right),
\left( \begin{array}{ccc}
1 & 0 & y \\
0 & 1 & 0 \\
0 & 0 & 1 \end{array} \right) \Bigg)
\mapsto
\left( \begin{array}{ccc}
1 & 0 & xy \\
0 & 1 & 0 \\
0 & 0 & 1 \end{array} \right)$$ is continuous with respect to $T'(GL_3(R_1),G)$. 
The latter follows, since
$$
\left( \begin{array}{ccc}
1 & 0 & xy \\
0 & 1 & 0 \\
0 & 0 & 1 \end{array} \right)=
\left( \begin{array}{ccc}
1 & x & 0 \\
0 & 1 & 0 \\
0 & 0 & 1 \end{array} \right)
\left( \begin{array}{ccc}
1 & 0 & 0 \\
0 & 1 & y \\
0 & 0 & 1 \end{array} \right)
\left( \begin{array}{ccc}
1 & -x & 0 \\
0 & 1 & 0 \\
0 & 0 & 1 \end{array} \right)
\left( \begin{array}{ccc}
1 & 0 & 0 \\
0 & 1 & -y \\
0 & 0 & 1 \end{array} \right),$$\\
and maps 
$$
\left( \begin{array}{ccc}
1 & 0 & x \\
0 & 1 & 0 \\
0 & 0 & 1 \end{array} \right)\mapsto
\left( \begin{array}{ccc}
1 & x & 0 \\
0 & 1 & 0 \\
0 & 0 & 1 \end{array} \right)=$$\\
$$=
\left( \begin{array}{ccc}
1 & 0 & x \\
0 & 1 & 0 \\
0 & 0 & 1 \end{array} \right)
\left( \begin{array}{ccc}
-1 & 0 & 0 \\
0 & 1 & 0 \\
0 & 1 & 1 \end{array} \right)
\left( \begin{array}{ccc}
1 & 0 & x \\
0 & 1 & 0 \\
0 & 0 & 1 \end{array} \right)
\left( \begin{array}{ccc}
-1 & 0 & 0 \\
0 & 1 & 0 \\
0 & -1 & 1 \end{array} \right),$$ \\
$$ 
\left( \begin{array}{ccc}
1 & 0 & y \\
0 & 1 & 0 \\
0 & 0 & 1 \end{array} \right) \mapsto
\left( \begin{array}{ccc}
1 & 0 & 0 \\
0 & 1 & y \\
0 & 0 & 1 \end{array} \right)=$$\\
$$=\left( \begin{array}{ccc}
1 & 0 & y \\
0 & 1 & 0 \\
0 & 0 & 1 \end{array} \right)
\left( \begin{array}{ccc}
1 & 0 & 0 \\
1 & 1 & 0 \\
0 & 0 & -1 \end{array} \right)
\left( \begin{array}{ccc}
1 & 0 & y \\
0 & 1 & 0 \\
0 & 0 & 1 \end{array} \right)
\left( \begin{array}{ccc}
1 & 0 & 0 \\
-1 & 1 & 0 \\
0 & 0 & -1 \end{array} \right),$$\\
$$\left( \begin{array}{ccc}
1 & 0 & z \\
0 & 1 & 0 \\
0 & 0 & 1 \end{array} \right)
\mapsto
\left( \begin{array}{ccc}
1 & 0 & -z\\
0 & 1 & 0 \\
0 & 0 & 1 \end{array} \right)$$

are continuous.

\hfill $\square$\\
The proof of the theorem has been completed.
\hfill $\square$\\

\begin{remark}
In the context of Theorem \ref{ring_top}, we obtain the same topology 
on $R$ if we identify $R$ with
$$\{\left( \begin{array}{ccc}
1 & x & 0 \\
0 & 1 & 0 \\
0 & 0 & 1 \end{array} \right):x\in R\}.$$ 
\end{remark}
{\em Proof.} This follows from the fact that 

$$\left( \begin{array}{ccc}
1 & x & 0 \\
0 & 1 & 0 \\
0 & 0 & 1 \end{array} \right)\mapsto
\left( \begin{array}{ccc}
1 & 0 & x \\
0 & 1 & 0 \\
0 & 0 & 1 \end{array} \right)=
\left( \begin{array}{ccc}
1 & -x & x \\
0 & -1 & 1 \\
0 & 0 & 1 \end{array} \right)
^2=\Bigg(
\left( \begin{array}{ccc}
1 & x & 0 \\
0 & 1 & 0 \\
0 & 0 & 1 \end{array} \right)
\left( \begin{array}{ccc}
1 & 0 & 0 \\
0 & -1 & 1 \\
0 & 0 & 1 \end{array} \right)
\Bigg)^2$$
is a continuous function (the continuity of the inverse to that
map follows from the calculations made at the end of the proof of 
Theorem \ref{ring_top}).
\hfill $\square$\\

\setcounter{claim}{0}
\section{Topologies on Polish structures}\label{topologies}

In this section, we will study the topologies introduced above
in the context of Polish structures,
which were introduced in \cite{31}, and studied also in \cite{35,82,86,34}.

\begin{definition}\label{smallness}
A Polish structure is a pair $(X,G)$, where $G$ is a Polish group acting faithfully on a set $X$ 
so that the stabilizers of all singletons are closed subgroups of $G$. If $X$ is equipped with
a structure of a group, which is preserved under the action of $G$, then we call $(X,G)$
a Polish group structure.
We say that $(X,G)$ is small if for every $n<\omega$, there are only countably 
many orbits on $X^n$ under the action of $G$.
\end{definition}
The class of small Polish structures contains examples of the form $(X,Homeo(X))$ (where
$Homeo(X)$ is considered with the compact-open topology) for $X$ being one of the spaces
$[0,1]^n,S^n$  ($n$-dimensional sphere), $(S^1)^n$ for $n\in\omega \cup \{\omega\}$,
as well as various other examples, see \cite [Chapter 4]{31}.

By Remark \ref{tau}, we have:
\begin{corollary}\label{taupol}
(1) If $G$ is a Polish group acting on a set $X$, then $(X,G)$ is a Polish structure iff
$\tau (X,G)$ is $T_1$ iff $\tau (X,G)$ is completely metrizable.\\
(2) If $(X,G)$ is a small Polish structure, then $\tau(X,G)$ is a Polish topology. In particular,
$(X,G)$ is a Polish $G$-space if we equip $X$ with the topology $\tau(X,G)$.
\end{corollary}
Let $(X,G)$ be a Polish structure. For any finite $C\subseteq X$, by $G_C$ we denote the pointwise stabilizer of $C$ in $G$, and for a finite tuple $a$ of elements of $X$, by $o(a/C)$ we denote the orbit of $a$ under the action of $G_C$ (and we call it the orbit of $a$ over $C$).

A fundamental concept for \cite{31} is the relation of $nm$-independence in an arbitrary Polish structure.
\begin{definition}
 Let $a$ be a finite tuple and $A$, $B$ finite subsets of $X$. Let $\pi _A:G_A\to o(a/A)$ be defined by $\pi _A(g)=ga$. We say that $a$ is $nm$-independent from $B$ over $A$ (written $a\nmind_A B$) if $\pi _A^{-1}[o(a/AB)]$ is non-meager in $\pi _A^{-1}[o(a/A)]$. Otherwise, we say that $a$ is $nm$-dependent on $B$ over $A$ (written $a\nmdep_A B$).
\end{definition}
By \cite[Theorem 2.14]{31}, under some assumptions, $nm$-dependence in a $G$-group $(H,G)$ can be expressed in terms of the topology on $H$:
\begin{fact}\label{gdelta}
Let $(X,G)$ be a Polish structure such that $G$ acts continuously on a Hausdorff space $X$. Let $a,A,B\subseteq X$ be finite. Assume that $o(a/A)$ is non-meager in its relative topology. Then, $a\nmind_A B \iff o(a/AB)\subseteq _{nm}o(a/A)$.
\end{fact}

Using the above fact, we now express the relation of $nm$-independece in terms of a 
familiy of topologies on $X$, without assuming anything about the Polish structure $(X,G)$.
\begin{remark}\label{express}
Let $(X,G)$ be a Polish structure and let $a,A,B \subseteq X$ be finite. 
Then $a\nmind_{A}B \iff o(a/B) \subseteq_{nm} o(a/A)$, where $X$ is equipped with the
topology $\tau(G_A,X)$ (and
the action of $G_A$ on $X$ is the restriction of the action of $G$ on $X$).
\end{remark}
{\em Proof.} The conclusion follows from Fact \ref{gdelta} and Corollary \ref{taupol}(2).
\hfill $\square$\\

If $A$ is a finite subset of $X$ (where $(X,G)$ is a Polish structure), we define the algebraic closure of $A$ (written $Acl(A)$) as the set of all elements of $X$ with countable orbits over $A$. If $A$ is infinite, we define $Acl(A)=\bigcup\{Acl(A_0):A_0\subseteq A$ is finite$\}$.
By Theorems 2.5 and 2.10 from \cite{31}, $nm$-independece has some nice properties 
corresponding to those of forking independece in stable first-order theories:
\begin{fact}\label{properties}
In any Polish structure $(X,G)$, $nm$-independence has the following properties:\\
(0) (Invariance) $a\nmind_A B\iff g(a)\nmind_{g[A]} g[B]$ whenever $g\in G$ and $a,A,B\subseteq X$ are finite.\\
(1) (Symmetry) $a\nmind_C b\iff b\nmind_C a$ for every finite $a,b,C\subseteq X$.\\
(2) (Transitivity) $a\nmind_B C$ and $a\nmind_A B$ iff $a\nmind_A C$ for every finite $A\subseteq B\subseteq C\subseteq X$ and $a\subseteq X$.\\
(3) For every finite $A\subseteq X$, $a\in Acl(A)$ iff for all finite $B\subseteq X$ we have $a\nmind_A B$.\\
If additionally $(X,G)$ is small, then we also have:\\
(4) (Existence of $nm$-independent extensions) For all finite $a\subseteq X$ and $A\subseteq B\subseteq X$ there is $b\in o(a/A)$ such that $b\nmind_A B$.
\end{fact}

Using Remark \ref{express}, we can slightly simplify some of the arguments from \cite{31}. 
For example, we reprove the existence of non-forking extensions in small Polish 
structures (point 4 of Fact \ref{properties}):\\
Let $a\subseteq X$ 
and $A\subseteq B\subseteq X$ be all finite. Since $\tau(X,G_A)$ is Polish, and 
there are countably many orbits over $B$, we can find, by the Baire category theorem, an element $b\in o(a/A)$, 
such that $o(b/B)$ is non-meager in  $\tau(X,G_A)$. Then, by Remark \ref{express}, 
we get that $b\nmind_{A}B$.

We will now apply Corollary \ref{description} to some of the structures 
constructed in \cite[Chapter 2]{82}. First, we outline the construction of those structures.

Suppose $(X,G)$ is a Polish structure. Let $H$ be an arbitrary group. 
For any $x\in X$ we consider an isomorphic copy $H_x=\{h_x:h\in H\}$ of $H$. 
By $H(X)$ we will denote the group $\bigoplus_{x\in X} H_x$. 
Although $H(X)$ is not necessarily commutative, we will denote its group action by $+$. 
For any $y\in H(X)$ there are $h_1,\dots,h_n\in H\backslash\{e\}$ and pairwise 
distinct $x_1,\dots,x_n\in X$ such that $y=(h_1)_{x_1}+\dots +(h_n)_{x_n}$. 
We will then write $h(y)=\{x_i:h_i=h\}$.

The group $G$ acts as automorphisms on $H(X)$ by $$g((h_1)_{x_1}+\dots +(h_n)_{x_n})=(h_1)_{gx_1}+\dots +(h_n)_{gx_n}.$$ 
It was proved in \cite{82} that with this action $(H(X),G)$ is a Polish structure, 
and that if $H$ is countable, and $(X,G)$ is small, then
also $(H(X),G))$ is small. Moreover, it was proved there  that if, additionally, $X$
is uncoutable, then these structures do not possess
any $nm$-generic orbits (the notion of an $nm$-generic orbit was 
introduced in \cite[Definition 5.3]{31}).  
On the other hand, \cite[Theorem 5.5]{31} states:

\begin{fact}\label{generics}
Suppose $(H,G)$ is a small Polish group structure, where $H$ is equipped with a topology
in which $H$ is not meager in itself (and the action of $G$ on $H$ is continuous).
Then, at least one $nm$-generic orbit in $H$ exists, and an orbit is $nm$-generic in $H$ iff 
it is non-meager.
\end{fact}

From the above theorem and the absence of generics it was concluded that for any non-trivial
countable group $H$
and any small Polish structure $(X,G)$, if $X$ is uncountable then there is no
non-meager in itself Hausdorff group topology on $H(X)$, such that the action
of $G$ on $H(X)$ is continuous (in particular, there is no such Polish topology). 
We strengthen this observation in some cases:

\begin{proposition}\label{haus}
Let $H$ be any non-trivial group and let $X$ be a compact Hausdorff
space containing an open subset homeomorphic to $(0,1)^n$ for 
some non-zero $n\in \omega\cup\{ \omega\}$ (notice that the
 examples listed after
Definition \ref{smallness} satisfy this assumption). Then, there is no 
Hausdorff group topology on 
$H(X)$ under which the action of $Homeo(X)$ on $H(X)$ is continuous
(where $Homeo(X)$ is considered with the compact-open topology).
\end{proposition}
{\em Proof.} Suppose first that $X=[0,1]$.
It is enough to show that the topology $\rho:=T(H([0,1]),Homeo([0,1]))$ is not Hausdorff. 
Take any $a\in H\backslash \{e\}$. We will show that any $\rho$-open neighbourhood of $e\in H([0,1])$
contains the element $a_{1/3}-a_{2/3}$. Let $W$ be any such neighbourhood and choose (by 
Corollary \ref{description}) a $\rho$-open set 
$$V=\bigcup_{q_1<\dots<q_n\in\mathbb{Q}}\{h_1(u_1h_1^{-1})u_1(h_2(u_2h_2^{-1}))u_1u_2(h_3(u_3h_3^{-1}))\dots u_1u_2\dots u_{n-1}(h_n(u_nh_n^{-1})):$$\\ 
$$h_i\in H,u_i\in U_{h_i}^{q_i}, u_1\dots u_n=e\},$$ such that $V+V\subseteq W$.
Let $B_{\epsilon}(id)\subseteq Homeo([0,1])$ be a ball (in the supremum metric) contained in 
$U^0_0\cap U^3_0$, and choose $n<\omega$ such that $1/3n<\epsilon$. 
Put $$h= a_{n/3n}+a_{(n+1)/3n}+
\dots +a_{(2n-1)/3n},h'= -a_{(n+1)/3n}-a_{(n+2)/3n}-\dots -a_{2n/3n}$$ and $U=U^1_{h}\cap
U^1_{h'}\cap U^2_0$. Notice that 
$\{u_0(h-u_1h), u_0(h'-u_1h'):u_0\in B_{\epsilon}(id),
u_1\in U\}\subseteq V$ (to see this, choose $q_j=j$ for $j=0,1,2,3$,
$u_2=u_1^{-1}$, $u_3=u_0^{-1}$, $h_0=h_2=h_3=0$ and $h_1$ equal to either $h$ or
$h'$).
Since $U$ is open, we can find $u_1\in U$ such that $u_1(k/3n)\in 
(k/3n,(k+1)/3n)$ for $k=n,n+1,\dots, 2n-1$. Then, there is some $u_0\in B_{\epsilon}(id)$
such that $u_0(k/3n)=k/3n$ and $u_0u_1(k/3n)=(2k+1)/6n$ for $k=n,n+1,\dots, 2n-1$. So,
we get that $\sum_{k=n}^{2n-1}(a_{k/3n}-a_{(2k+1)/6n})\in V$. 
Similarly, we obtain using $h'$
that $\sum_{k=n+1}^{2n}(-a_{k/3n}+a_{(2k-1)/6n})\in V$.  
Thus, $a_{1/3}-a_{2/3}\in V+V\subseteq W$.

Now, suppose $X$ is any space as in the statement. Then, we can find a copy $F$ of $[0,1]^n$ 
contained in $(0,1)^n\subseteq X$, and an isometric (with respect to a fixed metric on $F$)
copy $I$ of $[0,1]$ contained in $F$, such that
every homeomorphism of $I$ preserving its endpoints can be extended to a
homeomorphism of $F$ having the 
same distance from the identity (with respect to the supremum metric) and equal to the 
identity on the border of $F$ in $(0,1)^n$. 
Futhermore, since $(0,1)^n$ is open in $X$,
we can extend such a homeomorphism of $F$ to a homeomorphism of $X$ equal
to the identity on $X\backslash F$. Notice that any open neighbourhood 
of $id\in Homeo(X)$
contains $\{f\in Homeo(X):f_{|X\backslash int(F)}=id,d(id_F, f_{|F})<\epsilon\}$ for 
some $\epsilon>0$, where $d$ is the supremum metric.
Indeed, by the definition of the compact-open topology, such a neighbourhood is 
of the form $\{f\in Homeo(X):f[K_1]\subseteq W_1,\dots, f[K_l]\subseteq W_l\}$
where $W_i$'s are open, and each $K_i$ is a compact subset of $W_i$.
Then, it is enough to choose $\epsilon$ such that for each $i$ and $x\in K_i\cap F$,
$B_F(x,\epsilon)\subseteq F\cap U_i$.
Now, we can repeat the proof that we gave in the case of $X=[0,1]$.
Namely, choosing $V$ as above (but for an arbitrary $X$), we define $\epsilon$ to be 
such that $\{f\in Homeo(X):f_{|X\backslash int(F)}=id,d(id_F, f_{|F})<\epsilon\}
\subseteq U^0_0\cap U^3_0$ and define $h,h'$ in the same way as above 
(identyfing $[0,1]$ with a subset of $F$). Since $u_0$ and $u_1$
can be chosen to preserve endpoints of $[0,1]$, the choice of $F$ and of
the copy of $[0,1]$ inside it allows us to repeat the argument.
\hfill $\square$\\

The only known examples of small Polish group structures without $nm$-generic orbits
are of the form $H(X)$. For those of them for which we were able to 
compute the finest compatible topology, it turned out that it is not Hausdorff.
This may suggest that there could be a topological property of 
a group $H$ other than being non-meager in itself, which guarantees
the existence of $nm$-generic orbits in a structure $(H,G)$.

\begin{problem}
Characterize the existence of $nm$-generic orbits in a Polish group structure
$(H,G)$ in terms of topological properties of $H$.
\end{problem}

In particular, we can ask:

\begin{question}
Does the existence of a Hausdorff group topology on a group $H$ such that the action
of a Polish group $G$ on $H$ is continuous imply that the structure $(H,G)$ 
has an $nm$-generic orbit?
\end{question}

Also, we do not know whether the converse is true.

\begin{question}
Does the existence of $nm$-generic orbits in a Polish group structure
$(H,G)$ imply the existence of
a compatible Hausdorff topology on $H$?
\end{question}


\noindent
{\bf Address:}\\
Instytut Matematyczny, Uniwersytet Wroc\l awski,\\
pl. Grunwaldzki 2/4, 50-384 Wroc\l aw, Poland.\\[3mm]
{\bf E-mail address:}\\
dobrowol@math.uni.wroc.pl \\


\begin{thebibliography}{9999}



\bibitem{83} G. M. Bergman, 
{\em 
On group topologies determined by families of sets}, preprint, 2013.


\bibitem{35} R. Camerlo,
{\em Dendrites as Polish structures},Proceedings of the American Mathematical Society
(139), 2217-2225, 2011.


\bibitem{82} J. Dobrowolski,
{\em New examples of small Polish structures}, Journal of Symbolic Logic (78), 969-976, 2013.

\bibitem{86} J. Dobrowolski, K. Krupi\'nski,
{\em Locally finite profinite rings}, Journal of Algebra (401), 161-178, 2014. 

\bibitem{87} A. S. Kechris,
{\em Classical Descriptive Set Theory}, Springer, New York, 1995.




\bibitem{31} K. Krupi\'nski,
{\em Some model theory of Polish structures}, 
Transactions of the American Mathematica Society (362), 3499-3533, 2010.

\bibitem{34} K. Krupi\'nski, F. Wagner,
{\em Small, nm-stable compact G-groups}, 
Israel Journal of Mathematics (194), 907-933, 2013. 



\end{thebibliography}
\end{document}